\documentstyle[leqno,12pt,twoside]{article}

\makeatletter
\makeatletter

\@addtoreset{equation}{section}

\newcommand{\decl}{:=}                  
\newcommand{\abs}[1]{\left|#1\right|}   
\newcommand{\be}{\begin{equation}}      
\newcommand{\ee}{\end{equation}}        
\newcommand{\bern}{\begin{eqnarray*}}   
\newcommand{\eern}{\end{eqnarray*}}     
\newcommand{\uni}{{\bf 1}}

\newtheorem{defi}{Definition}[section]
\newtheorem{theorem}[defi]{Theorem}
\newtheorem{prop}[defi]{Proposition}
\newtheorem{lemma}[defi]{Lemma}
\newtheorem{rem}[defi]{Remark}
\newtheorem{coro}[defi]{Corollary}

\def\bbbr{{\rm I\!R}} 

\def\bbbn{{\rm I\!N}} 

\def\L{L_{n,\lambda_n}}

\begin{document}

\thispagestyle{empty}
\par\vskip 35mm

\begin{center} {\large SOME REMARKS ON A GENERAL CONSTRUCTION OF  APPROXIMATION PROCESSES} \end{center}

\begin{center}Lorenzo D'Ambrosio,  Elisabetta  Mangino 

Ref. SISSA 8/2001/M
\end{center}

\par\vskip  15mm

The  aim of this paper is to study the characteristics of a general method
to produce  a new approximation sequence from a given one, by using suitable
convex combinations, whose coefficients depend on  some functions $\lambda_n$
($n\ge 1$).
This type of construction was first used to generate Lototsky-Schnabl
operators from Bernstein operators (see section 6.1. in \cite {AC}) and
turned to be really useful in the study of evolution
equations by means of positive operators.
More precisely,  the new sequences of operators were used
to represent the solutions of degenerate  elliptic parabolic equations,
thus providing several qualitative informations both on the solutions
of the equation and on the Markov process associated with the equation
(see  Chapter 6 in \cite{AC} and the references quoted therein and e.g. \cite{Ac,Am1,Am2}).

In the present note we investigate which properties of the original
sequence of operators
$(L_n)_{n\ge 1}$ are inherited by the new sequence $(\L)_{n\ge 1}$.
More precisely, approximation properties and rates of approximations are studied in
the first section. The second section is devoted to the study
of  the asymptotic behaviour of the  sequence $(\L)_{n\ge 1}$ and a
Voronovskaja-type relation is established.
Finally in the third section the behaviour of the operators $\L$, when acting on convex,
monotone or H\"older continuous functions is  studied.

\bigskip

Throughout the paper, $I$ will denote an interval on the real line and $J$ a subinterval of $I$.
Let
$F(I)$ be the space of all real valued functions on $I$ and
let $C(I)$ denote the space of all
real valued continuous functions on $I$. Moreover  set
 $C_b(I)$ for be the subspace of all  bounded
continuous functions on $I$.

The previous spaces of bounded functions,
endowed with the natural order and the norm
\[ \Vert f\Vert_\infty:=\sup_{x\in I}\vert f(x)\vert, \qquad f\in
 C_b(I),\]
are  Banach lattices.

If $g:I\rightarrow \bbbr$ is a strictly positive function, let $C(I,g)$ be the subspace
of $C(I)$ of all functions $f$ such that $f/g$ is bounded.

Set $e_n(x):=x^n$ for every $x\in I$ and $n\in\bbbn$, and
denote by ${\bf 1}$
the function with constant value $ 1$. Moreover, for every $x\in I$,
let $\psi_x(t):=(t-x)$ ($t\in I$).
For every $f\in C(I)$, for every
$\alpha\in [0,1]$ and $x\in I$, set
\[ f_{\alpha,x}(t):=f\left( \alpha t + \left(1-\alpha\right)x\right)
\qquad (t\in I).\]

Let  $D$ be a  sublattice of $C(I)$ such
that  ${\bf 1},e_1\in D$ and    for every
$f\in D$, $x\in I$, $\alpha\in [0,1]$, the function $f_{\alpha,x}$ belongs
to $D$.

In this paper $(L_n)_{n\ge 1}$ will denote a sequence of linear positive operators from $D$
into $F(J)$ such that  $L_n({\bf 1})={\bf 1}$ ($n\ge 1$). We set
$L_0(f):=f_{|J}$ for every $f\in D$.

For every $n\ge 1$, we consider $\lambda_n:J\rightarrow [0,1]$,
and  for every $f\in D$ and $x\in J$ we define the operator
\[ L_{n,\lambda_n}(f)(x):=\sum_{p=0}^n {n\choose p} \lambda_n(x)^p
(1-\lambda_n(x))^{n-p}L_p(f_{p/n,x})(x). \]

Clearly  $L_{n,\lambda_n}$ is a
linear positive operator from $D$ into  $F(J)$.

If $\lambda_n(x)=1$ (resp. $\lambda_n(x)=0$) for some $x\in J$,
then $L_{n,\lambda_n}(f)(x)=L_n(f)(x)$   (resp.
$L_{n,\lambda_n}(f)(x)=f(x)$).

\section{Approximation properties}

In this section we investigate under which conditions the sequence
$(L_{n,\lambda_n})_{n\ge 1}$ is an approximation process.
For this we will apply some Korovkin-type theorems and their extension to
spaces of unbounded functions (see   \cite{AC, LD, ditzian, shaw} and
references therein).

The first step is to evaluate the operators $L_{n,\lambda_n}$ on
some ``test'' functions.

\begin{lemma}\label{l1} For every $n\ge 1$,
        \be \L({\bf 1})={\bf 1}.\label{t0}\ee

        If $\psi_x^k \in D$, for $x\in J$ and $k\ge 1$, then
        \be \L(\psi_x^k)(x)=\frac1n \sum_{p=0}^{n}{{n}\choose p}
        \lambda_n(x)^p(1-\lambda_n(x))^{n-p}\frac{p^k}{n^{k-1}}
        L_{p}(\psi_x^k)(x). \label{tk}\ee
        In particular, for any $x\in J$
        \be \L(\psi_x)(x)=\lambda_n(x)\sum_{p=0}^{n-1}{{n-1}\choose p}
        \lambda_n(x)^p(1-\lambda_n(x))^{n-1-p}L_{p+1}(\psi_x)(x),\label{t1}\ee
        and,  if $e_2\in D$,
        \begin{eqnarray}\label{t2}
& & \L(\psi_x^2)(x)= \\
&=& \frac{\lambda_n(x)}{n}
\sum_{p=0}^{n-1}
        {{n-1}\choose p} \lambda_n(x)^p(1-\lambda_n(x))^{n-1-p}(p+1)
                L_{p+1}(\psi_x^2)(x). \nonumber \end{eqnarray}
\end{lemma}
{\bf Proof.} The assertions follow by easy computations, once we observe
that  \begin{eqnarray*} {\uni}_{\frac p n,x}&=&\uni,\nonumber\\
        (\psi_x^k)_{\frac p n,x} &=& \frac{p^k}{n^k} \psi_x^k, \nonumber
\end{eqnarray*}
for any $x\in I$, $0\le p\le n$ and $k\ge 0$.
\hfill$\Box$

\medskip
\begin{rem} If $S:D\to F(J)$ is a positive linear operator and $e_2\in D$,
        then for any $x\in J$, we have
        \be \abs{S(\psi_x)(x)}\le (1+S(\uni)(x)) \sqrt{S(\psi_x^2)(x)}.
                                        \label{sat}\ee
        Indeed,  the relation
        \be \abs{\psi_x}\le (\uni+\psi_x^2\delta^{-2})\delta,
\label{eas}\ee
 holds        for any $\delta>0$. By applying $S$ to (\ref{eas}) and then
choosing         $\delta:=S(\psi_x^2)(x)$, we obtain the claim.
\end{rem}

In the sequel we will consider some assumptions on the operators $L_n$.
Let $(S_n)_{n\ge 1}$ be a sequence of positive linear operator from $D$
into $F(J)$ and $M_2:J\to\bbbr_+$ be a positive function.
We say that $(S_n)_{n\ge 1}$ satisfies the condition $({\mathbf
H}_2)$ with growth  $M_2$, if
\begin{enumerate}
\item[$({\mathbf H}_2)$] $e_2\in D$  and
$S_n(\psi_x^2)(x)\le \frac{M_2(x)}{n}$ for any $n\ge 1$ and $x\in J$.
\end{enumerate}
If $({\mathbf H}_2)$ is satisfied, then,  by (\ref{sat}), for every
$x\in J$ we have \be\sup_{n\in\bbbn}\abs{\sqrt n S_n(\psi_x)(x)}\le
                (1+S_n(\uni)(x))\sqrt{M_2(x)}.\label{hh}\ee

Sometimes we shall  need of a stronger condition than (\ref{hh}).
Namely, let $M_1:J\to\bbbr_+$ be a positive function,
the sequence $(S_n)_{n\ge 1}$ satisfies the condition $({\mathbf
H}_1)$ with growth $M_1$, if
\begin{enumerate}
\item[$({\mathbf H}_1)$] $\abs{S_n(\psi_x)(x)} \le
\frac{M_1(x)}{n}$         for any $n\ge 1$ and $x\in J$.
\end{enumerate}

For instance, the conditions $({\mathbf H}_i)$ are satisfied whenever a
Voronovskaja-type relation holds (see next section).

Notice that, if for any $n\ge 1$ $L_n$ preserves the affine functions,
then from (\ref{t0}) and (\ref{t1}) we deduce that the same
property holds for any $\L$. In this case $({\mathbf H}_1)$ is verified
with  $M_1=0$.

The conditions $({\mathbf H}_i)$
are preserved by the transformation
$(L_n)_{n\ge 1} \rightarrow (\L)_{n\ge 1}$, as specified by the following
\begin{lemma}
\begin{enumerate}\item  If $(L_n)_{n\ge 1}$ satisfies $({\mathbf H}_2)$
        with growth $M_2$, then for any $n\ge 1$ and $x\in J$ we have
        \be \L(\psi_x^2)(x)\le \frac{\lambda_n(x)M_2(x)}{n},\label{1}\ee
        therefore, also $(\L)_{n\ge1}$ satisfies $({\mathbf H}_2)$ with
	growth $M_2$,  and moreover
        \be\abs{\L(\psi_x)(x)}\le 2\sqrt{\frac{\lambda_n(x) M_2(x)}{n}}.
                        \label{x}\ee

\item If $(L_n)_{n\ge 1}$ satisfies $({\mathbf H}_1)$ with growth $M_1$,
then         for any $n\ge 1$ and $x\in J$ we have
        \be\abs{\L(\psi_x)(x)}\le M_1(x)\frac{1-(1-\lambda_n(x))^n}{n},
                                \label{2}\ee
        therefore, also $(\L)_{n\ge1}$ satisfies $({\mathbf H}_1)$ with
growth $M_1$. \end{enumerate}
\end{lemma}
{\bf Proof.} {\em 1.}  The hypothesis $({\mathbf H}_{2})$ applied to $L_p$
in         (\ref{t2}) of Lemma \ref{l1} yields (\ref{1}).
        Combining (\ref{sat}) and (\ref{1}), we obtain (\ref{x}).

{\em 2.} By (\ref{t1}),  we obtain
\[\abs{\L(\psi_x)(x)}\le \lambda_n(x) \sum_{p=0}^{n-1} {n-1 \choose p}
        \lambda_n(x)^p (1-\lambda_n(x))^{n-1-p}\frac{M_1(x)}{p+1}.\]
Since
\[ s \sum_{p=0}^{n-1} {n-1\choose p}\frac{ s^p (1-s)^{n-1-p}}{p+1} =
        \frac{1}{n} \sum_{p=0}^{n-1} {n \choose p+1}s^{p+1} (1-s)^{n-1-p}=
        \frac{1-(1-s)^{n}}{n}, \]
holds for  $0\le s \le 1$, we obtain (\ref{2}) by replacing
$s$ with $\lambda_n(x)$.
\hfill$\Box$\medskip

\begin{prop} Assume that $C_b(I)\subseteq D$ and $(L_n)_{n\ge 1}$ satisfies
        $({\mathbf H}_2)$ with growth $M_2$.
\begin{enumerate}
\item For every $f\in C_b(I)$, $n\ge1$, $x\in J$ the following estimate holds:
        \[ \vert\L(f)(x)-f(x)\vert \leq 2\omega\left(f, \sqrt{{
        {\lambda_n(x)M_2(x)}\over n}}\right).\]

\item For every differentiable $f\in D$ such that $f'\in C_b(I)$,
        $n\ge1$, $x\in J$ the following estimate holds:
        \[ \vert L_{n,\lambda_n} (f)(x)-f(x)\vert \leq
        2\left(\abs{f'(x)} +
        \omega(f',  \sqrt{{{\lambda_n(x)M_2(x)}\over n}})\right)
        \sqrt{{{\lambda_n(x)M_2(x)}\over n}},\]
	and if $({\mathbf H}_1)$ is satisfied with growth $M_1$,
        then we have
\begin{eqnarray*} \vert L_{n,\lambda_n} (f)(x)-f(x)\vert &\leq&
        \vert f'(x)\vert M_1(x)\frac{ 1-(1-\lambda_n(x))^n}{n} +\\
        &&\qquad + 2\omega(f',  \sqrt{{{\lambda_n(x)M_2(x)}\over n}})
        \sqrt{{{\lambda_n(x)M_2(x)}\over n}}.\end{eqnarray*}

\item If $I=J=(0,\infty)$, $(L_n)_{n\ge 1}$ satisfies
        $({\mathbf H}_1)$ with growth $M_1$ bounded and $({\mathbf H}_2)$
with growth  $M_2(x)=\alpha +\beta x+\gamma x^2$ ($x\in I$),
 $\alpha,\beta,\gamma\in\bbbr$, and moreover we assume
        $\L(C_b(I))\subset C_b(I)$, then for every $f\in C_b(I)$,
        $n\ge1$, we have
        \[  \Vert L_{n,\lambda_n}(f)-f\Vert_\infty \leq
        K\left(\omega_{\sqrt{M_2}}^2(f,{{1}\over \sqrt n})
                +{{1}\over n}\right), \]
where $K>0$ is a constant independent of $n$, and
\[ \omega_\phi(f,\delta):=\sup_{0\leq h\leq \delta, x\pm
h \phi(x)\geq 0} \vert \Delta^2_{h\phi(x)}f(x)\vert\]
with
\[ \Delta^2_{h\phi(x)}f(x)= f(x-h\phi(x))-2f(x)+f(x+h\phi(x)).\]
If $L_n(e_1)=e_1$, the previous estimate becomes
\[  \Vert L_{n,\lambda_n}(f)-f\Vert_\infty \leq
        \omega_{\sqrt{M_2}}^2(f,\sqrt{{\Vert \lambda_n\Vert}\over n}).\]
\end{enumerate}
\end{prop}
{\bf Proof.}
The statements  in {\em 1.} and {\em 2.} follow from
Theorem
5.1.2 and the subsequent remark in \cite{AC}, taking into account  the
estimates (\ref{1}), (\ref{x}) and (\ref{2}).

{\em 3.} is an immediate consequence of Theorem 1 and the subsequent remark
in \cite{t} and the estimates (\ref{1}) and (\ref{2}). \hfill$\Box$
\medskip
\begin{rem} Note that in the setting of 3. of previous proposition,
if $\Vert\lambda_n\Vert\rightarrow 0$, there also holds
\[  \Vert L_{n,\lambda_n}(f)-f\Vert_\infty \leq
        K\left(\omega_{\sqrt{M_2}}^2(f,\sqrt{\Vert\lambda_n\Vert})
                +\Vert\lambda_n\Vert\right). \]
Indeed using Bernoulli inequality $(1+s)^n\ge 1+ns$ ($s\ge-1$), from (\ref{2})
one has $\abs{\L(\psi_x)(x)}\le M_1(x)\frac{1-(1-\lambda_n(x))^n}{n}\le
M_1(x)\lambda_n(x)$, and then, we can argue as before.
\end{rem}

In the sequel, with $g:I\to \bbbr$ we shall denote a strictly convex
function such that $g\ge c$  for some constant $c>0$.

\begin{rem}\label{oss} \begin{enumerate}
\item If $f$ is (strictly) convex, then also $f_{\alpha,x}$ is (strictly) convex,
        for any $\alpha\in [0,1]$ and $x\in I$.
\item If $f\in C(I,g)$, then also $f_{\alpha,x}\in C(I,g)$,
        for any $\alpha\in [0,1]$ and $x\in I$.
\item If $f\in C^k(I)$, then also $f_{\alpha,x}\in C^k(I)$,
        for any $\alpha\in [0,1]$ and $x\in I$, $k=0,1,2,\dots$
\end{enumerate}
\end{rem}

\begin{theorem}\label{conv_growth}  Let $g:I\to \bbbr$ be a
strictly convex function unbounded on $I$
with a continuous derivative  $g'\colon J\to \bbbr$
such that $g\ge c$  for some constant $c>0$.
        In case $I$ is unbounded, we additionally require
        \[ \lim_{{\scriptstyle \abs t\rightarrow\infty}}
                \frac{g(t)}{\abs{t}}=+\infty.\label{eqn:ipo.cresc.inf.g}\]
We assume that $C(I,g)\subset D$ and the following conditions
$({\mathbf G})$ and $({\mathbf H}_g)$ hold,
\begin{enumerate}
\item[$({\mathbf G})$] $L_p(g_{\frac p n,x})(x)\ge g_{\frac p n,x}(x)=g(x)$
        for any $x\in J$, $n\ge 1$ and $p=0,\dots ,n$,
\item[$({\mathbf H}_g)$]
        $M_g(x):=\sup_{p\in\bbbn}p\abs{L_p(g)(x)-g(x)}$
        is finite for every $x\in J$.

\item If $(L_n)_{n\ge 1}$ satisfies $({\mathbf H}_1)$ with growth $M_1$,
        then for all $f\in C(I,g)$, we have
\be \lim_{n\to\infty}L_{n,\lambda_n}(f)=f \label{con}\ee
pointwise, and  the convergence is uniform on the sets where $M_1$,
$M_g$ and $g'$ are bounded.

\item  If $(L_n)_{n\ge 1}$ satisfies $({\mathbf H}_2)$ with growth $M_2$,
        then there exists a constant $M>0$,
        depending only on $f$, $I$ and $J$, such that for every $n\ge1$,
        $x\in J$, we have
\begin{eqnarray}
\vert \L(f)(x)-f(x)\vert &\leq& 2\omega\left(f,
                \sqrt{{{\lambda_n(x)M_2(x)}\over n}}\right) +\nonumber \\
        &&\qquad +M\left[{{\lambda_n(x)M_g(x)}\over n}+2 \abs{g'(x)}
        \sqrt{\frac{ M_2(x)\lambda_n(x)}{n}}\right]. \label{1dis}
\end{eqnarray}
        Moreover, if also $({\mathbf H}_1)$ holds with growth $M_1$,
        then for $n\ge1$, $x\in J$, we  have
\begin{eqnarray}
\vert \L(f)(x)-f(x)\vert&\leq& 2\omega\left(f,
        \sqrt{{{\lambda_n(x)M_2(x)}\over n}}\right) + \nonumber \\
        && +M {{\lambda_n(x)M_g(x)}+\abs{g'(x)}M_1(x) (1-(1-\lambda_n(x))^n)  \over n}. \label{2dis}
\end{eqnarray}
\end{enumerate}
\end{theorem}
{\bf Proof.} We start with some consequence of conditions
$({\mathbf G})$ and $({\mathbf H}_g)$.

The identity
\be \qquad \L(f)(x)-f(x)=\sum_{p=0}^n {n\choose p} \lambda_n(x)^p
        (1-\lambda_n(x))^{n-p}\left(L_p(f_{p/n,x})(x)-f(x)\right),
                \label{id}\ee
and $({\mathbf G})$ assures that $\L(g)(x)-g(x)\ge 0$.

The convexity of $g$  yields
\[ g_{\frac p n,x}(t)-g(x)\le \frac p n (g(t)-g(x)),\]
which with identity (\ref{id}) and $({\mathbf H}_{g})$ implies
\bern 0\le \L(g)(x)-g(x) &\le& \sum_{p=1}^n {n\choose p} \lambda_n(x)^p
                (1-\lambda_n(x))^{n-p}\frac p n(L_p(g)(x)-g(x))\\
        &\le& \frac{M_g(x)}{n}\sum_{p=1}^n {n\choose p} \lambda_n(x)^p
                (1-\lambda_n(x))^{n-p}= \frac{M_g(x)\lambda_n(x)}{n}.
\eern

{\em 1.} Relation (\ref{2}) together with the previous one,
        assure the convergence of $\L$ on the ``test''
        functions $\uni$, $e_1$ and $g$.
        Now, using Theorem 4.1 in \cite{LD}, we obtain the convergence (\ref{con}).

{\em 2.} In this case the relation (\ref{x}) and (\ref{1}) hold.
        Thus using Theorem 4.5 in \cite{LD} we get (\ref{1dis}).
        Moreover, if even $({\mathbf H}_{1})$ holds, then (\ref{2})
and again Theorem 4.5 in \cite{LD}, yield (\ref{2dis}).
\hfill$\Box$

\medskip

Observe that  condition $({\mathbf G})$ is satisfied whenever $L_n$ has
the following property: \[ f {\mathrm \ is\ convex\ }\Rightarrow  L_n(f)\ge
f.\] Another concrete condition implying $({\mathbf G})$ is provided by the
following lemma.

\begin{lemma}\label{lem_conv} We assume that $L_n$ preserves the affine function, that is
         $L_n(e_i)=e_i$ for $i=0,1$. For any convex function
        $f\in C^2(I)\cap D$, we have $L_n(f)\ge f$.
\end{lemma}
{\bf Proof.} Indeed, expanding $f$ in Taylor series, we have
\[ f(u)=f(x)+f'(x)(u-x)+\frac{f''(\eta)}{2}(u-x)^2\ge f(x)+f'(x)(u-x),\]
and applying $L_n$ to both side of the inequality, we have $L_n(f)\ge f$.
\hfill $\Box$

Actually, also  condition $({\mathbf H}_g)$ is easily checkable.
For instance, it is satisfied if a Voronovskaja-type relation holds for
$L_n$ (see next example).

\medskip
\noindent{\bf Example.}
Set $EXP(I)\decl\bigcup_{w>0} C(I,\exp(we_1)+\exp(-we_1))$,
the space of continuous functions with exponentially growth.
For any $n\ge 1$ let
$L_n\colon EXP(I)\to C(J)$ be a linear positive operator defined by
\[ L_n(f)(x)=\int_I W_n(x,t)f(t) d t,\]
where the kernel $W_n(x,t)$ is a generalized positive function.
We assume that there exists a strictly positive function $p\in C(I)$,
analytic on the inner of $I$, such that the relations
\begin{eqnarray} L_n(\uni)(x)&=&1,\nonumber\\
                \frac{\partial}{\partial x}L_n(f)(x)&=&\frac{n}{p(x)}
                L_n(\psi_xf)(x),\label{eqn:exp.oper.deriv}
\end{eqnarray}
hold for any $f\in EXP(I)$ and $x\in J$.
These operators are also referred as the \emph{exponential operators}
(for more details see \cite{ismail,may}).
Operators satisfying these property are, for example,
the Bernstein polynomials, the Sz\'sz-Mirakjan, Baskakov,
Post-Widder and Weierstrass operators.

These operators realize approximation processes for functions in
$EXP(I)$, and the following results hold (see \cite{ismail, may}).
\begin{prop}\label{prop_exp} For any $n\ge 1$, $x\in J$ and
$f\in EXP(I)\cap C^2(I)$, the following assertions hold:
\begin{enumerate}
\item $L_n(\uni)(x)=1$,
\item $L_n(e_1)(x)=x$,
\item $L_n(e_2)(x)=x^2+\frac{p(x)}{n}$,
\item $\lim_n n(L_n(f)(x)-f(x))=\frac 12 p(x)f''(x)$ uniformly on any
        compact subset of the interior of $J$.
\end{enumerate}
\end{prop}

Now, we can state the following theorem.
\begin{theorem} For any $f\in C(I)$, such that
$\abs{f(x)}\le K(\exp(wx)+\exp(-wx))$,
there exists a constant $M(f)>0$, depending only on $f$, for which
the estimate
\[ \abs{\L(f)(x)-f(x)}\leq 2\omega\left(f, \sqrt{{{\lambda_n(x)p(x)}\over
n}}\right) + M(f)w^2 {{\lambda_n(x)}p(x)\cosh(wx)\over n},\]
holds uniformly on compact subset of $J$.
\end{theorem}
{\bf Proof.} By Proposition \ref{prop_exp},
        $(L_n)_{n\ge1}$ satisfies $({\mathbf H}_2)$ with $M_2=p$ and
        $({\mathbf H}_1)$ with $M_1=0$.
        We set $g(x)=\exp(wx)+\exp(-wx)$.
        From Lemma \ref{lem_conv} and 4. of Proposition \ref{prop_exp},
        we deduce respectively that $({\mathbf G})$ and $({\mathbf H}_g)$
        are  satisfied.
        Thus, applying Theorem \ref{conv_growth}, we obtain the claim.
\hfill $\Box$

\section{Voronovskaja-type results}

\begin{prop}\label{vo}
Assume $J$ compact, there exist
$\alpha,\beta \in C(J)$, an even number $q>2$ such that $e_k\in D$ for
$k=2,\dots, q$, and   for every $x\in J$, there hold

(i)  $\lim_{n\to\infty} nL_n(\psi_x)(x)=\beta(x)$,

(ii) $\lim_{n\to\infty} nL_n(\psi_x^2)(x)=\alpha(x)$,

(iii) $ \lim_{n\to\infty} nL_n(\psi_x^q)(x)=0$.

Moreover assume that the sequence
$(\lambda_n)_{n\ge 1}$ converges uniformly to a function $\lambda$, and
$\lambda_n(x)>0$ for every $n\ge 1$  and $x\in J$. Then
for every $f\in C^2(I)\cap D$, such that $f''$ is bounded,  the convergence
\begin{equation}\label{vo1}
\lim_{n\to\infty} n(\L(f)(x)-f(x))=\lambda(x)
{{\alpha(x)}\over 2}f''(x)+\beta(x)f'(x) \end{equation}
holds pointwise on $J$.

In particular, if $\lambda=0$, then
\[ \lim_{n\to\infty} n(\L(f)(x)-f(x))=\beta(x)f'(x). \]
\end{prop}
{\bf Proof.}
Let $f\in C^2(I)\cap D$ bounded on $I$ and $x\in J$. Set
$M:=\max_{I}f''$.Then

$$f=f(x){\bf 1} + f'(x)\psi_x + {{f''(x)}\over 2}\psi_x^2 +
\omega(x,\cdot)\psi_x^2,$$
where $\omega(x,y)={{f''(\xi)-f''(x)}\over 2}$ with
$\xi$ between $x$ and $y$.
Hence
\[ n(\L(f)(x)-f(x))= f'(x)n\L(\psi_x)(x)+{{f''(x)}\over 2}n\L(\psi_x^2) +
n\L(\omega(x,\cdot)\psi_x^2)(x). \]

By \ref{tk}  in Lemma \ref{l1} and
by Toeplitz theorem (see e.g. \cite{st}, Theorem (7.85)), as $n$ tends to
infinity, $n\L(\psi_x)(x)$ tends to
$\beta (x)$ and $n\L(\psi_x^2)(x)$ tends to $\alpha(x)$.
Moreover $f''$ is  continuous in $x$, thus, for a fixed
$\varepsilon>0$ there exists $\delta >0$ such that $\vert x-y\vert<\delta$
implies  $\vert f''(x)-f''(y)\vert<\varepsilon$, for every $y\in I$.

Therefore, if $\vert x-y\vert<\delta$, then
$\vert \omega(x,y)\psi_x^2(y)\vert\leq {\varepsilon}\psi_x^2(y)$, while, if
$\vert x-y\vert\geq \delta$, then
$\vert \omega(x,y)\psi_x^2(y)\vert \leq M\delta^{2-q}\psi_x^q(y)$.
Hence
\[ \vert n\L(\omega(x,\cdot)\psi_x^2)(x)\vert \leq \varepsilon
n\L(\psi_x^2)(x)+M\delta^{2-q}n\L(\psi_x^q)(x).\]
Thus, by applying again Toeplitz theorem, we get
$\lim_{n\to\infty}n\L(\omega(x,\cdot)\psi_x^2)(x)=0$. \hfill $\Box$

\bigskip

\noindent {\bf Example 1.} The convergence in the Voronovskaja-formula
        (\ref{vo1}) is uniform  on $J$ when
        $\lambda_n(x)\geq \lambda_0$ for every $n\ge 1$,
        $x\in J$, and the convergence in (i), (ii), (iii)
        are uniform on $J$, since in this case
        Toeplitz theorem yields the uniform
        convergence (see examples after Theorem (7.85) in \cite{st}).
        But the convergence can be uniform even when $\lambda=0$.

        For this we discuss the case of Bernstein-Schurer operators
        $B_{n,1}:C([0,2])\rightarrow C([0,1])$, defined for every
        $f\in C([0,2])$ and $x\in [0,1]$ by
\[ B_{n,1}(f):=\sum_{k=0}^{n+1}{{n+1}\choose k} x^k(1-x)^{n+1-k}f
        \left({k\over n}\right) \]
(see e.g. 5.3.1 in \cite{AC}).

It is known that $B_{n,1}(\psi_x)(x)={1\over n}x$, and
$B_{n,1}(\psi_x^2)(x)={1\over n}x(1-x)+{x\over{n^2}}$.
Simple calculations show that for every
$x\in [0,1]$, $B_{n,1}(\psi_x^4)(x)\leq
{M\over {n^3}}$, where $M$ is a constant that does not depend on $n$ nor on
$x$.

Choose $\lambda_n={1\over n}$ for every $n\ge 1$. Then
\[ B_{n,1,\lambda_n}(\psi_x)(x)= {x\over n}, \quad
B_{n,1,\lambda_n}(\psi_x^2)(x) = {x\over n^2} + {1\over n^2}x(1-x),\quad
B_{n,1,\lambda_n}(\psi_x^4)(x)\leq {M\over n^3}. \]
By reviewing the proof of \ref{vo}, one concludes that for every
$f\in C^2([0,2])$,
\[ \lim_{n\to\infty}n(B_{n,1,\lambda_n}(f)-f)(x)=xf'(x)\]
uniformly on $[0,1]$.

\bigskip
\noindent{\bf Example 2.} Let $(L_n)_{n\ge 1}$ be the sequence
of exponential operators that have been considered in the example in
        the previous section. By combining (\ref{eqn:exp.oper.deriv})
        and  {\em 3.} of Proposition \ref{prop_exp}, one has
        $\lim_n nL_n(\psi_x^4)(x)=0$.
        Now, in the setting of Proposition \ref{vo}, that is,
        assuming $J$ compact, $\lambda_n>0$ for any $n\ge1$ and
        $\lambda_n\rightarrow \lambda$ uniformly,
        we have that, for every $f\in C^2(I)\cap EXP(I)$ such that
        $f''\in C_b(I)$, there holds
        \[ \lim_n n(\L(f)(x)-f(x))=\frac{\lambda(x)p(x)}{2}f''(x) \quad
(x\in J).\]

\section{Shape preserving properties}

We investigate the behaviour of the sequence  $(L_{n,\lambda_n})_{n\ge 1}$
when acting on convex, monotone or H\"older continuous functions. The
results
do not depend  on the topological structure of the space $D$, but only on
 the positivity of the operators $L_n$
and their shape preserving properties.

\begin{prop}  For every $n\ge 1$, let  $\lambda_n=\lambda: J\to[0,1]$.
Assume that for every convex function $f\in D$, the sequence $(L_n(f))_{n\ge 1}$ is decreasing.
Then $(L_{n,\lambda}(f))_{n\ge 1}$ is also decreasing for every
convex function $f\in D$.
\end{prop}
{\bf Proof.} Let $f\in D$ be a convex function and $x\in J$.
        We have to prove that the quantity
        $Q\decl L_{n+1,\lambda}(f)(x)-L_{n,\lambda}(f)(x)=
        L_{n+1,\lambda}(f)(x)-L_{n,\lambda}(f)(x)(1-\lambda(x)+\lambda(x))$
        is not positive.
        For sake of shortness, in the sequel of this proof we avoid to
        indicate the dependence from $x$  of $\lambda$, $f$ and $L_p$.
        By computation one gets
\begin{eqnarray*}
Q\!&\!=\!&\!\sum_{p=0}^{n+1} {{n+1}\choose p} \lambda^p(1-\lambda)^{n+1-p}L_p(f_{p/(n+1),x})-
 \sum_{p=0}^{n} {{n}\choose p} \lambda^{p+1}(1-\lambda)^{n-p}L_p(f_{p/n,x})+\\
&& -\sum_{p=0}^{n} {n\choose p} \lambda^p(1-\lambda)^{n+1-p}L_p(f_{p/n,x})=\\
\!&\!=\!&\! (1-\lambda)^{n+1}f+ \lambda^{n+1} L_{n+1}(f)-(1-\lambda)^{n+1}f-\lambda^{n+1}L_n(f) +\\
&&+\left[ (n+1)L_1(f_{1/(n+1),x})-f-nL_1(f_{1/n,x})\right]\lambda(1-\lambda)^n
        +\sum_{p=2}^n \lambda^p(1-\lambda)^{n+1-p}\times\\
&&\times\left[ {{n+1}\choose p} L_p(f_{p/(n+1),x})-
{n\choose {p-1}} L_{p-1}(f_{(p-1)/n,x})- {n\choose p} L_p(f_{p/n,x})\right].
\end{eqnarray*}

For every $n\ge 1$ and $p=0,...,n-1$ the function $f_{p/n,x}$
is convex, therefore $L_p(f_{p/n,x})\geq L_{p+1}(f_{p/n, x})$.
By the positivity of the operators $L_p$ and the convexity of $f$, we obtain
\begin{eqnarray*}
Q\!&\!\leq\!&\!\lambda(1-\lambda)^n\left[ (n+1)\left({n\over{n+1}}L_1(f_{1/n,x})+
        {1\over {n+1}} f \right) -f-nL_1(f_{1/n,x})\right]+\\
\!&+\!&\!\sum_{p=2}^n\left[ {{n+1}\choose p} L_p(f_{\frac{p}{n+1},x})-
        {n\choose {p-1}} L_{p}(f_{\frac{p-1}{n},x})- {n\choose p}
        L_p(f_{\frac{p}{n},x})\right] \lambda^p(1-\lambda)^{n+1-p}\\
\!&\!=\!&\!\sum_{p=2}^nL_p\left[ {{n+1}\choose p} f_{\frac{p}{n+1},x}-
        {n\choose {p-1}} f_{\frac{p-1}{n},x}- {n\choose p}
        f_{\frac{p}{n},x}\right]\lambda^p(1-\lambda)^{n+1-p}.
\end{eqnarray*}
From the convexity of $f$ and positivity of $L_p$, we get $Q\leq0$.
\hfill$\Box$

\begin{coro}  If for every convex function $f\in D$ and $n\ge 1$,
\[ f\leq L_{n+1}(f)\leq L_n(f),\]
then for  every convex function $f\in D$ and $n\ge 1$, we have
\[ f\leq L_{n+1,\lambda}(f)\leq L_{n,\lambda}(f)\leq \lambda L_1(f) +
        ({\bf 1}-\lambda)(f) \leq L_1(f).\]
\end{coro}
{\bf Proof.} Let $f\in D$ be convex, then $f_{p/n,x}$ is convex for every
        $n\ge 1$ and $p=0,...,n$. Hence for every $t\in J$
        \[ L_p(f_{p/n,x})(t)\geq f_{p/n,x}(t).\]
        In particular, $L_p(f_{p/n,x})(x)\geq f(x)$, hence
        $L_{n,\lambda}(f)(x)\geq \sum_{p=0}^n{n\choose p}\lambda(x)^p(1-\lambda(x))^{n-p}f(x)=f(x)$.
        The other inequalities follow by observing that
        $L_{1,\lambda}=\lambda  L_1+ ({\bf 1}-\lambda)Id_D$
        and that $f\leq L_1(f)$. \hfill $\Box$

\bigskip

Consider the following subsets of $F(I)$:
\begin{eqnarray*}
&{\bf M}_+(I):=\{ f:I\rightarrow \bbbr: f\  {\rm is\ increasing} \},\\
&{\bf Cx}(I):=\{f:I\rightarrow \bbbr:  f\ {\rm is\ convex}\},\\
\!\!&\!\!{\bf Lip}_K(\alpha)(I):= \{ f:I\rightarrow \bbbr:  \vert f(x)-f(y)\vert
        \leq K\vert x-y\vert^\alpha\ (x,y\in I)\} \ (K>0, \alpha>0).
\end{eqnarray*}

\begin{prop} Let $\lambda_n \in [0,1]$ be a constant.
\begin{enumerate}
\item If $L_p({\bf  M}_+(I)\cap D))\subseteq {\bf  M}_+(J)$ for every
$p=1,...,n$, then $\L({\bf  M}_+(I)\cap D)\subseteq {\bf  M}_+(J)$.
\item If  $L_p({\bf  M}_+(I)\cap D)\subseteq {\bf  M}_+(J)$ and
        $L_p({\bf Cx}(I)\cap D)\subseteq {\bf Cx}(J)$ for every $p=1,...,n$,
        then $\L({\bf Cx}(I)\cap D)\subseteq {\bf Cx}(J)$.

\item If  $L_p({\bf Lip}_K(\alpha)(I)\cap D)\subseteq {\bf Lip}_H(\alpha)(J)$
        for every $p=1,...,n$, then
        $\L({\bf Lip}_K(\alpha)(I)\cap D)\subseteq {\bf Lip}_{K+H}(\alpha)(J)$.
        In particular Lipschitz  continuous functions are preserved.
\end{enumerate}
\end{prop}
{\bf Proof.} {\em 1.} See the proof of Proposition 2.2 in \cite{Ca}.

{\em 2.}  Let $f\in {\bf Cx}(I)\cap D$ and $x,y\in I$, $x< y$.
        As already stated in Remark \ref{oss}  $f_{p/n,x}\in {\bf Cx}(I)$
        for every $p=0,...,n$ and $x\in I$.
        By following  the idea of the proof of Proposition 2.3
        in \cite{Ca}, one gets
        \[ L_p(f_{p/n,y})(y)+L_p(f_{p/n,y})(x)\leq
                L_p(f_{p/n,x})(x)+L_p(f_{p/n,y})(y).\]

        Now, since  the operator $L_p$ is positive and
 preserves the convexity, we get for every $t\in[0,1]$,
\begin{eqnarray*}%
\lefteqn{L_p(f_{p/n,tx+(1-t)y})(tx+(1-t)y)\le}\\
	&&\qquad\qquad \le tL_p(f_{p/n,tx+(1-t)y})(x) +(1-t)L_p(f_{p/n,tx+(1-t)y})(y)\\
        &&\qquad\qquad \leq t^2 L_p(f_{p/n,x})(x) + t(1-t)L_p(f_{p/n,y})(y) +\\
        &&\qquad\qquad\qquad\qquad +(1-t)tL_p(f_{p/n,x})(x) +(1-t)^2L_p(f_{p/n,y})(y)\\
        &&\qquad\qquad =tL_p(f_{p/n,x})(x)+ (1-t)L_p(f_{p/n,y})(y).
\end{eqnarray*}
        Thus, the function $z\rightarrow L_p(f_{p/n,z})(z)$ is convex.
        It follows immediately that  $L_{n,\lambda_n}(f)$ is convex too.

{\em 3.} Let $f\in {\bf Lip}_K(\alpha)(I)\cap D$ and $x,y\in I$.
        Then $f_{p/n,x}\in {\bf Lip}_{K(p/n)^\alpha}(\alpha)(I)$ and hence
\begin{eqnarray*}
        \lefteqn{\vert \L(f)(x)- L_{n,\lambda_n}(f)(y)\vert \leq
        	\sum_{p=0}^n{n\choose p}\lambda_n^p(1-\lambda_n)^{n-p}\times}\\
        &&\qquad\times\left[ \vert L_p(f_{p/n,x})(x)-L_p(f_{p/n,x})(y)\vert +
        \vert L_p(f_{p/n,y})(y)-L_p(f_{p/n,x})(y)\vert\right] \\
        &&\leq \sum_{p=0}^n{n\choose p} \lambda_n^p(1-\lambda_n)^{n-p}
        \left[ K\left({p\over n}\right)^\alpha +
        H L_p( \left( 1-{p\over n}\right)^\alpha \vert x-y\vert^\alpha \cdot
        {\bf 1}) (y)\right]\\
        &&\leq (H+K)\vert x-y\vert^\alpha,
\end{eqnarray*}
that is the claim. \hfill$\Box$

\begin{rem} If $\lambda_n$ is not a constant, the previous results do not hold in general.
In \cite{Ca}, there are several counterexamples in the case of
Sz\'asz-Mirakjan operators. \end{rem}

Finally, we do some remarks about the monotonicity of the operators
$L_{n,\lambda_n}$ with respect to the function $\lambda_n$.
If $\lambda\le \alpha$,  in general we do not have
$L_{n,\lambda}(f)\leq L_{n,\alpha}(f)$, even
when $\alpha$ and $\lambda$ are constants and $f$ is convex.
For a counterexample consider $\lambda={1\over 2}$, $\alpha=1$,
$L_n$ be the n-th Kantorovich operator on $[0,1]$ and $f=e_1$ (see
\cite{AC}, 5.3.7).

\begin{prop} Let $\alpha_n, \lambda_n:J\rightarrow [0,1]$  be functions
        such that $\alpha_n\leq \lambda_n$.
        If $f\in D$ and $x\in J$ satisfy
        $L_p(f_{p/n,x})(x)\leq L_{p+1}(f_{(p+1)/n,x})(x)$
        for every $p=0,...,n-1$,
        then $L_{n,\alpha_n}(f)(x)\leq L_{n,\lambda_n}(f)(x)$.
\end{prop}

{\bf Proof.} Let $f\in D$ and $x\in J$
satisfying the assumptions. The function
\[ \phi:\ s\in [0,1] \rightarrow \sum_{p=0}^n {n\choose p} s^p(1-s)^{n-p} L_p(f_{p/n,x})(x), \]
has  derivative
\[ \phi'(s)=n\sum_{p=0}^n {{n-1}\choose p} s^p(1-s)^{n-p-1} (L_{p+1}(f_{(p+1)/n,x})(x)- L_p(f_{p/n,x})(x)) \quad (s\in [0,1]),\]
hence it is increasing. In particular $L_{n,\alpha_n}(f)(x)=\phi(\alpha_n(x))\leq \phi(\lambda_n(x))=L_{n,\lambda_n}(f)(x)$.\hfill $\Box$

\begin{rem} If for every $x\in J$, there exists a regular
         finite Borel measure $\gamma_x$ on $I$ such that
        \[ L_p(f)(x)= \int_{I^p} f\left( {{x_1+\dots+x_p}\over p}
        \right)d\gamma_x(x_1)\cdots d\gamma_x(x_p) \qquad (f\in D),\]
        and $f_{|J}\leq L_1(f)$ for every convex function $f$, then
        $L_p(f_{p/n,x})(x)\leq L_{p+1}(f_{(p+1)/n,x})(x)$ for every convex
        function $f$, $x\in J$, $n\ge 1$ and $p=1,...,n-1$.

Indeed, let $x_0,x\in J$, then for every $\alpha\in [0,1]$, we have
        \[ f\left(\alpha x_0+(1-\alpha)x\right)= f_{\alpha,x_0}(x) \leq
        L_1(f_{\alpha,x})(x)=\int_I f\left(\alpha x_0 + (1-\alpha)t\right)
        d\gamma_x(t).\]
        Since $\displaystyle{{{x_1+\dots+x_p +x}\over{p+1}} }\in I$,
        whenever $x_1,\dots,x_p,x\in I$, we get
\begin{eqnarray*}
        \lefteqn{ L_p(f_{p/n,x})(x)= \int_{I^p} f\left( {{x_1+...+x_p+x}\over
        {p+1}}{{p+1}\over n} + {{(n-p-1)x}\over n}\right) d\gamma_x(x_1)
                \cdots d\gamma_x(x_p)}\\
        &&\leq\! \int_{I^p}\!\left(\int_I f\left({{x_1+...+x_p+x_{p+1}}\over n}+
        {{(n-p-1)x}\over n}\right)\! d\gamma_x(x_{p+1})\right) d\gamma_x(x_1)
        \cdots d\gamma_x(x_p) \\
        &&= L_{p+1}(f_{(p+1)/n,x})(x).
\end{eqnarray*}

Particular cases of this situation are Bernstein-Schnabl operators (see
\cite{AC}), Szasz-Mirakjan and  Baskakov operators on $[0,+\infty[$ and
their modification of Lototskj-Schnabl  type (see \cite{AC,Ac,Am1}).
\end{rem}

\bigskip

Author's address:

LORENZO D'AMBROSIO

SISSA - ISAS v. Beirut, 2-4

I-34014  TRIESTE  ITALY

E-MAIL:dambros@sissa.it

\bigskip

ELISABETTA  MANGINO

DIPARTIMENTO DI MATEMATICA ``E. DE GIORGI''

UNIVERSITA' DEGLI STUDI DI LECCE

I-73100 LECCE

E-MAIL:elisabetta.mangino@unile.it
\end{document}